\documentclass{amsart}

\usepackage{amsthm,amsfonts,amsmath,amssymb,mathrsfs,verbatim}
\usepackage{hyperref}

\numberwithin{equation}{section}
\newtheorem{theorem}{Theorem}[section]
\newtheorem{corollary}[theorem]{Corollary}
\newtheorem{proposition}[theorem]{Proposition}

\newtheorem{question}[theorem]{Question}

\newtheorem{lemma}[theorem]{Lemma}
\newtheorem{remark}[theorem]{Remark}
\newtheorem{definition}[theorem]{Definition}
\newtheorem{example}[theorem]{Example}

\newcommand{\Z}{\mathbb Z}
\newcommand{\F}{\mathbb F}

\DeclareMathOperator{\GCD}{GCD}

\DeclareMathOperator{\PF}{PF}

\begin{document}

\title{Star operation on orders in simple Artinian rings}

\author{Nazer H. Halimi}

\thanks{}

\email{n.halimi@uq.edu.au}

\address{School of Mathematics and Physics, 
The University of Queensland, QLD 4072, Australia}

\date{August 2011}

\begin{abstract}
Star operations are an important tool in multiplicative ideal theory. 
In this paper we apply a special type of star operation,
known as $\nu$-operation, to define the notion of
right Pr\"ufer $\nu$-multiplication order.
The latter may be viewed as a natural non-commutative version of 
Pr\"ufer $\nu$-multiplication domain.
As one of our main results, we establish that  
an overring of a right Pr\"ufer $\nu$-multiplication 
order is again a right Pr\"ufer $\nu$-multiplication order.
\end{abstract}

\maketitle 

\section{Introduction} 
Multiplicative ideal theory is a crucial ingredient in the 
classification of orders in simple Artinian rings \cite{IMU04,M09,MMU97}.
In turn, star operations are a powerful tool used to study
multiplicative ideal theory. Most progress, however, 
is concerned with the application of star operations
in the commutative setting, see 
\cite{AAMZ09,EB,FZ09,PT} and references therein,
and relatively little is known in the non-commutative case.
To the best of our knowledge, the first to advance the latter
were Asano and Murada \cite{AM53}, who, in 1953, used the 
$\nu$-operation to study the Arithmetic properties of 
non-commutative semigroups. 
Surprisingly, only a handful of further works have been devoted to
the application of star operations in non-commutative ring theory.
For example, in \cite{BD03,BMO00} the $\nu$-operation is used
to classify prime segments of Dubrovin valuation rings, and right cones 
in right orderable groups.
Furthermore, in \cite{M09}, Marubayashi applied the $\nu$-operation 
to investigate Ore extensions over total valuation rings. 
He defined the notion of $\nu$-Bezout orders, and proved that an order in a
simple Artinian ring is $\nu$-Bezout if and only if it
is a $\GCD$ order.

In this article we build on Marubayashi's work, to further
our understanding of $\nu$-multiplication orders.
We study both the $\nu$-Bezout order introduced in \cite{M09}
as well as a new order, called 
right Pr\"ufer $\nu$-multiplication order, which is a
non-commutative version of Pr\"ufer $\nu$-multiplication domain.
For this new order we prove that an overring of a right Pr\"ufer 
$\nu$-multiplication 
order is again a right Pr\"ufer $\nu$-multiplication order.

The remainder of the paper organised as follows.
In the next section our main theme of study is the $\nu$-Bezout order. 
In particular, in Theorem \ref{total}, we show that for 
total valuation rings $V\subset W$ of a division ring $K$, 
the $\nu$-Bezout order $V+W[x,\sigma]x$ in $K(x,\sigma)$ is Bezout if 
and only if $W=K$.
In Section~\ref{sec3}, our focus is the $\tau$-operation, 
which, among other things, is needed in our subsequent study
of right Pr\"ufer $\nu$-multiplication order in Section~\ref{sec4}.
%However, we generalized some results from commutative to non-commutative case
%by applying star operation in general. 
In that particular section, we first define this new order
and then prove that an overring of a right Pr\"ufer $\nu$-multiplication 
is itself a right Pr\"ufer $\nu$-multiplication order.      
%Throughout this paper ring $S$ will be a prime Goldie ring having total quotient ring $Q$.

\section{$\nu$-Bezout orders}\label{sec2}

Let $S$ be an order in a simple Artinian ring $Q$, i.e., $S$ is a prime Goldie ring with total quotient ring $Q$.
Given regular elements $a, b$ in $S$, we say that $b$ is a \textit{right-divisor} of $a$,
if there exists an element $s\in S$ such that $a=bs$. A \textit{left-divisor} is defined similarly.
\begin{definition}
Let $a, b, d\in S$ be regular elements of $S$. We say that $d$ is a \textit{right greatest common divisor}
of $a$ and $b$, denoted by $d=r-\gcd\{a, b\}$, if the following two conditions hold:
\begin{itemize}
\item[(i)]
$d$ is a right-divisor of $a$ and $b$;
\item[(ii)]
if $c$ is a right-divisor of $a$ and $b$, then $c$ is a right-divisor of $d$.
\end{itemize}
\end{definition} 
Note that if $c$ is another right greatest common divisor of $a$ and $b$, then $dS=cS$.

A left greatest common divisor of $a, b$, denoted by $l-{\gcd}\{a, b\}$, is defined likewise.
An order $S$ in a simple Artinian ring $Q$ is called a $\GCD$ order if any two elements of $S$ have a right 
as well as a left greatest common divisor.

Let $U(Q)$ is the group of units in $Q$. A right $S$-submodule $I$ of $Q$ is called a \textit{right $S$-ideal}
if $I$ contains a regular element in $S$ and $uI\subseteq S$ for some $u\in U(Q)$. For any subsets
$A$ and $B$ of $Q$, we use the notations:
$(A:B)_r=\{q\in Q: qB\subseteq A\}$ and $(A:B)_l=\{q\in Q:Bq\subseteq A\}$.
If $I$ is a right $S$-ideal, then $(S:I)_l$ is a left $S$-ideal. 
We define $I_{\nu}:=(S:(S:I)_l)_r$. The set $I_{\nu}$ is a right $S$-ideal containing $I$. 
If $I_{\nu}=I$, then it is called a \textit{right $\nu$-ideal}. Similarly, for any left $S$-ideal $J$, we can 
define a left $S$-ideal $_{\nu}J$. An order $S$ in $Q$ is called a \textit{right $\nu$-Bezout} order
if $I_{\nu}$ is right principal for any finitely generated right integral $S$-ideal $I$.
A \textit{left $\nu$-Bezout} order is defined similarly.
By a \textit{$\nu$-Bezout} order we mean it is a right as well as left $\nu$-Bezout order.
In \cite[Proposition 2.2] {M09} Marubayashi proved that an order $S$ in a simple Artinian ring 
$Q$ is a $\nu$-Bezout order if and only if $S$ is a $\GCD$ order.

\begin{lemma}\label{I}
Let $V$ be a total valuation ring of a division ring $K$ and 
let $I$ be an ideal of $T:=V+K[x,\sigma]x$. Then the following three statements are equivalent:
\begin{itemize}
\item[(1)]
$I\cap V\neq 0$;
\item[(2)] $K[x,\sigma]x\subset I$;
\item[(3)]
$IK[x,\sigma]=K[x,\sigma]$.
\end{itemize}
If any of the above conditions hold, then $I=(I\cap V)+K[x,\sigma]x=(I\cap V)T$.
\end{lemma}
\begin{proof}
The proof is similar to the commutative case, see \cite[Lemma 4.11]{CMZ78}.
\end{proof}
Let $V$ be a proper total valuation ring of a division ring $K$ and $\sigma$ be an automorphism of $V$.
Let $V[x, \sigma]$ be the ring of Ore extension over $V$ with multiplication $xa=\sigma(a)x$ for all $a\in V$. 
The automorphism $\sigma$ thus naturally extends to an automorphism of $K$. Then $K[x,\sigma]$ is a principal ideal
ring and has the division ring $K(x,\sigma)$ as quotient ring. In \cite{M09} it is proved 
that $V[x,\sigma]$ is a $\nu$-Bezout order in $K(x,\sigma)$, which is not 
Bezout. The following is a similar type of 
result but with a different proof strategy.

\begin{theorem}\label{total}
Let $V$ be a total valuation ring of a division ring $K$ and let $W$ be a proper overring of $V$.
Then the following three statements are equivalent.
\begin{itemize}
\item[(1)]
$V+W[x,\sigma]x$ is a $\nu$-Bezout order in $K(x,\sigma)$;

\item[(2)]
$W=K$;

\item[(3)]
$V+W[x,\sigma]x$ is a Bezout order in $K(x,\sigma)$.
\end{itemize}
\end{theorem}

\begin{proof}
$\boldsymbol{(1)\Rightarrow (2).}$\;
Since $W$ is a proper overring of the total valuation ring $V$, we have $W=V_S$, where $S=V-J(W)$
and $V_S$ is the localization of $V$ at the Ore system $S$.
Therefore, $V+W[x,\sigma]x=V+V_S[x,\sigma]x$. Put $V^{(S)}:=V+V_S[x,\sigma]x$.
To prove that $W=K$, it is enough to show that $J(W)=0$.
Proceeding by contradiction, assume there exists a nonzero element  $a\in J(W)$.
Then $aV^{(S)}\subseteq sV^{(S)}$ for all $s\in S$.
Since $x=ss^{-1}x$, we also have $xV^{(S)}\subseteq sV^{(S)}$ for all $s\in S$.
Now let $d=l-\gcd\{a, x\}$. Then $a=df$ and $x=dg$ for some $f, g\in V^{(S)}$.
From $x=dg$ and the fact that $d$ is a constant, we conclude that $d\notin J(W)$, which implies that $d\in S$.
Therefore, $xV^{(S)}\subseteq d^nV^{(S)}$ and $aV^{(S)}\subseteq d^nV^{(S)}$ for all $n\geq 0$. Now if $d$ is a non-unit 
in $V$, then this is in contradiction with $\gcd\{d^{-1}a, d^{-1}x\}\in U(V^{(S)})$. If $d$ is a unit in $V$, then 
this is in contradiction with $V\neq W$. Therefore our assumption on the existence of $a$ is false, and $J(W)=0$. 
This concludes the proof that (1) implies $W=K$.

\noindent
$\boldsymbol{(2)\Rightarrow (3).}$\;
Let $T=V+K[x,\sigma]x$. First we will show that each right ideal of $T$ is of the form $f(x)FT=f(x)(F+K[x,\sigma]x)$.
To see this, let $I$ be a right ideal $T$. If $IK[x,\sigma]=K[x,\sigma]$, then by Lemma~\ref{I}, $I\cap V\neq 0$ and 
$I=I\cap V+K[x,\sigma]x=(I\cap V)T$. Thus it is enough to consider $f(x)=1$ and
$F=I\cap V$.

Since $K[x,\sigma]$ is principal ideal ring, if $IK[x,\sigma]\neq K[x,\sigma]$ then
$IK[x,\sigma]=f(x)K[x,\sigma]$ for some $f(x)\in K[x,\sigma]$.
Hence there exist $a_i\in I$ and $h_i\in K[x,\sigma]$ such that $f(x)=\sum_{i=0}^na_ih_i$.
Let $h_i=\sum_jq_{i,j}x^j$. Then $f(x)=\sum_{i=0}^na_iq_{i,0}+\sum_{i=0}^na_ih'_i$, where 
$q_{i,0}\in K$ and $h'_i\in K[x,\sigma]x$. Since $V$ is a total valuation ring of $K$, there exist 
$t\in V$ such that $q_{i,0}t\in V$ for all $1\leq i\leq n$. 
Hence $f(x)t=\sum_{i=0}^n(a_iq_{i,0}t)+\sum_{i=0}^n(a_ih'_it)$. Since
$a_i\in I, q_{i,0}t\in V$ and $h'_it\in K[x,\sigma]x$, we have $f(x)t\in I$. Let 
\begin{equation}\label{defF}
F=\{t\in V:f(x)t\in I\}.
\end{equation}
Then $F$ is a nonzero right $V$-submodule $K$. From $f(x)F\subseteq I$ and $FT=F+K[x,\sigma]x$, we have $I\supseteq f(x)FT
=f(x)(F+K[x,\sigma]x)$. Conversely, let $h(x)\in I$. Then 
$h(x)=f(x)(q_0+q_1x+\dots+q_mx^m)$, where $q_i\in K$. Put $h'(x)=f(x)(q_1x+\dots+q_mx^m)$. Then
$h'(x)\in f(x)K[x,\sigma]x$ and $h(x)=f(x)q_0+h'(x)$. Since $f(x)(q_1x+\dots+q_mx^m)\subseteq I$, we have
$f(x)q_0=h(x)-h'(x)\in I$. Thus $q_0\in F$ and $h(x)\in f(x)(F+K[x,\sigma]x)$ so that $I\subseteq f(x)(F+K[x,\sigma]x)$.
Therefore, $I\supseteq f(x)FT\supseteq f(x)(F+K[x,\sigma]x)\supseteq I$, which shows that $I=f(x)FT=f(x)(F+K[x,\sigma]x)$. 

Now let $I$ be a finitely generated right ideal of $T$. Then $F$ defined in \eqref{defF} is
a finitely generated right $V$-module such that $I=f(x)FT$. 
By construction $F=d_1V+\dots+ d_nV$, where $d_1,\dots, d_n\in V$.
Therefore, $F=dV$ for some $d\in V$ and $I=f(x)dVT=f(x)dT$, which shows that $T$ is a right Bezout ring. Similarly
one can prove that $T$ is a left Bezout ring.  

\noindent
$\boldsymbol{(3)\Rightarrow (1).}$\;
A Bezout order is always a $\nu$-Bezout order.
\end{proof}  

\begin{lemma}
Let $S$ be a $\nu$-Bezout order in a division ring.
For any triple of regular elements $a, b, c\in R$, if $\gcd\{a, b\}=1$ then $\gcd\{a, bc\}=\gcd\{a, c\}$.
\end{lemma}

\begin{proof}
Suppose $d=l-\gcd\{a, c\}$ and $a=a'd, c=c'd$. Then by \cite [Lemma 2.1]{M09}, $l-\gcd\{a', c'\}=1$.
From $l-\gcd\{a', c'\}=1$ and $l-\gcd\{a', b\}=1$ we conclude that $l-\gcd\{a', bc'\}=1$. Again by 
\cite [Lemma 2.1]{M09}, $l-\gcd\{a, bc\}=d$, as desired. 
Likewise one can proof that $r-\gcd\{a, bc\}=r-\gcd\{a, c\}$.
\end{proof}

\begin{corollary}\label{prime}
For any $a, b, c\in S$, if $\gcd\{a, b\}=1$ and $a$ is a divisor of $bc$, then $a$ is a divisor of $c$.
\end{corollary} 

\begin{definition}
Let $S$ be a $\nu$-Bezout order in a division ring $D$. 
A prime ideal $P$ of $S$ is a left $\PF$-prime ideal if $l-\gcd\{a, b\}\in P$ for any pair of regular elements $a, b\in P$.
\end{definition}
\begin{theorem}
Suppose that $P$ is a completely prime ideal of a $\nu$-Bezout order $S$ such that $S$ is localizable at $P$. 
Then $P$ is $\PF$-prime if and only if $S_P$ is a total valuation ring of $D$. 
\end{theorem}
\begin{proof}
Let $P$ be a left $\PF$-prime ideal of $S$ and $x=ab^{-1}\in D$, where $a, b\in S$ such that $l-\gcd\{a, b\}=1$. Since
$1\notin P$, we have $a\notin P$ or $b\notin P$. Thus $x=ab^{-1}\in S_P$ or $x^{-1}=ba^{-1}\in S_P$.

Conversely, suppose $P$ is not a left $\PF$-prime. Then there exist $a, b\in P-\{0\}$ such that $l-\gcd\{a, b\}\notin P$.
Let $d=l-\gcd\{a, b\}$ and $x=ad^{-1}, y=bd^{-1}$. Then $x, y\in P$ and, by \cite [Lemma 2.1]{M09}, $l-\gcd\{x, y\}=1$. 
We will show that neither $xy^{-1}$ nor $yx^{-1}\notin S_P$. If $xy^{-1}\in S_P$ then $xy^{-1}=ts^{-1}$ for some 
$t\in S$ and $s\in S-P$. Since $S-P$ is an Ore set, there exist $t_1\in S$ and $s_1\in S-P$ such that $s_1t=t_1s$.
From the above we can conclude that $xy^{-1}=ts^{-1}=s_1^{-1}t_1$ and $s_1x=t_1y$. Thus $y$ is a \textit{left-divisor}
of $s_1x$. Since $l-\gcd\{x, y\}=1$, by Corollary~\ref{prime}, $y$ is a \textit{left-divisor} $s_1$. 
Therefore, $s_1=ky$ for some $k\in S$ and $s_1\in P$, a contradiction. 
By a similar reasoning it follows that $yx^{-1}\notin S_P$. Thus $S_P$ is not a total valuation ring.   
\end{proof}

\begin{corollary}
Let $S$ be a $\nu$-Bezout order in a division ring $D$ such that $S$ is localizable at every completely prime ideal. Then:
\begin{itemize}
\item[(i)]
Every completely prime ideal contained in a $\PF$-prime ideal is again a $\PF$-prime ideal.
\item[(ii)]
The set of all completely prime ideals contained in a $\PF$-prime ideal is linearly ordered, and 
hence the set of all $\PF$-prime ideals forms a tree.
\end{itemize}
\end{corollary}

\section{star operation on orders in a simple Artinian rings}\label{sec3}

Throughout of the rest of the paper $S$ is an order in a simple Artinian ring $Q$ and  
 $F_r(S)$ ($\bar{F}_r(S)$) are the set of nonzero right $S$-ideals ($S$-submodules) of $Q$.

\begin{definition}\label{star}
A mapping $I\rightarrow I^*$ of $\bar{F}_r(S)$ into $\bar{F}_r(S)$ is called 
a \textit{semistar operation} on $S$ if the following three conditions hold for 
all $u\in U(Q)$ and $I, J\in \bar{F}_r(S)$:
\begin{itemize}

\item[(1)]$(uI)^*=uI^*$;
\item[(2)] if $I\subseteq J$ then $I^*\subseteq J^*$;

\item[(3)]
$I\subseteq I^*$ and $ (I^*)^*=I^*$.
\end{itemize}
\end{definition}

When $S^*=S$ the restriction $*$ to $F_r(S)$ also satisfies to the above tree conditions and is called a \textit{star operation} on $S$.
An element $I\in F_r(S)$ is called a \textit{right star ideal} if $I=I^*$.

\begin{example}\label{T}
\begin{itemize}
\item[(i)]
$i_d :\bar{F}_r(S)\rightarrow \bar{F}_r(S)$ defined by $I^{i_d}=I$, the identity map on $\bar{F}_r(S)$, is a 
semistar operation on $S$.
\item[(ii)]
The map $\nu: \bar{F}_r(S)\rightarrow \bar{F}_r(S)$ defined by $I^{\nu }:=(S:(S:I)_l)_r$ is a semistar operation.

\item[(iii)]
If $*$ is a semistar operation on $S$ then we can define the map $*_f: \bar{F}_r(S)\rightarrow {F}_r(S)$
by $I^{*_f}:=\cup F^*, F\in \bar{F}_r(S)$ where the union runs over all 
finitely generated $F\subseteq I$ such that $F\in \bar{F}_r(S)$.
It is easy to see that $*_f$ ,which is called a \textit{semistar operation of finite
type} associated to $*$, is indeed a semistar operation. In particular, the semistar operation of finite
type associated to $\nu$ is denoted by $\tau$.
\end{itemize}
\end{example}

\begin{lemma}\label{TL}
Let $*$ be a star operation on $S$ and $\{F_{\alpha}\}$ a family of elements of $F_r(S)$.
Then $(\sum_{\alpha}F_{\alpha})^*=(\sum_{\alpha}F_{\alpha}^*)^*$.
\end{lemma}
\begin{proof}
The proof is the same as in the commutative case, see \cite[Proposition 32.2] {G92}.  
\end{proof}

\begin{lemma}
Let $S$ be an order in a simple Artinian ring $Q$ and $A, B$ right $S$-ideals. Then:
\begin{itemize}

\item[(1)]
If $A^{\tau}=B^{\tau}$ then $A^{\nu}=B^{\nu}$.
\item[(2)]
If $A$ is a $\nu$-ideal, then $A$ is a $\tau$-ideal.
\item[(3)]
If $S$ satisfies the ascending chain condition on the set of left $\nu$-ideals, then every 
right $\tau$-ideal is a right $\nu$-ideal.
\item[(4)]
$S$ satisfies the ascending chain condition on the set of left $\nu$-ideals if 
and only if it satisfies this condition on the set of left $\tau$-ideals.
\end{itemize}
\end{lemma}

\begin{proof}
\begin{itemize}
\item[(1)]
This follows from $A\subseteq A^{\tau}\subseteq A^{\nu}$ and 
$(A^{\nu})^{\nu}=A^{\nu}$. 

\item[(2)]
Since $A$ is a $\nu$-ideal and $A^{\tau}\subseteq A^{\nu}$, we have $A^{\tau}=A$.
\item[(3)]
Let $A$ be a right ideal and $\{B_i\}$ the family of right ideals of $S$ 
contained in $A$. Then $(S:A)_l\subseteq (S:B_i)_l$ for all $i$. 
Since each $(S:B_i)_l$ is a left $\nu$-ideal and $S$ satisfies the 
ascending chain condition on the set of all left $\nu$-ideals, there exists 
a minimal element in $\{(S:B)_l\}$, say $(S:B_n)_l$. If $B_n^{\nu}\neq A$, 
there exists an element $b\in A-B^{\nu}$. Let $B=B_n+bS$.
 Then $B$ is a finitely generated right ideal and $B\subseteq A$ such that $(S:B)_l\subset (S:B_n)_l$
 This contradicts the minimality of $(S:B_n)_l$ so that $B^{\nu}_n=A$.
 Now let $A$ be a right $\tau$-ideal. By the above there exists a finitely generated right $S$-ideal
 $B$ such that $B\subseteq A$ and $B^{\nu}=A^{\nu}$. Since $B^{\tau}=B^{\nu}$ and 
 $A$ is right $\tau$-ideal, we conclude that $A\subseteq A^{\nu}\subseteq B^{\nu}=B^{\tau}\subseteq A^{\tau}=A$.
 Thus $A^{\nu}=A$.
 \item[(4)]
 This follows from (3) and (4). \qedhere
\end{itemize}
\end{proof}
\begin{lemma}\label{chain} 
Let $A_1\subset A_2\subset\cdots $ be right $\tau$-ideals. Then 
$A=\cup_{i\geq 1}A_i$ is a right $\tau$-ideal.
\end{lemma}
\begin{proof}
Let $B$ be a finitely generated right $S$-ideal with $B\subset A$. Then there exists $A_n$ such that 
$B\subseteq A_n$. Thus $B^{\nu}\subseteq A^{\tau}_n=A_n$, so that $B^{\nu}\subseteq A$.
\end{proof}

\begin{corollary}
If every right $\tau$-ideal is a finitely generated $S$-ideal, then $S$ satisfies the ascending chain 
condition on the set of all right $\tau$-ideals.
\end{corollary}
\begin{proof}
If $S$ does not satisfy the ascending chain condition on the set of all right $\tau$-ideals,
then there exists an infinite chain of right $\tau$-ideals, $A_1\subset A_2\subset \cdots$.
By Lemma~\ref{chain}, $A=\cup _{i\geq 1}A_i$ is a $\tau$-ideal which is not a finitely generated, a contradiction. 
\end{proof}
\begin{lemma}
Let $V$ be a total valuation ring of rank of least $2$ in a division ring $K$, and $Q$ a non-maximal completely
prime ideal of $V$. Put $T=V+V_Q[x,\sigma]x$. Then $M=\{f\in T: f(0)\notin U(V) \}$ is a $\tau$-ideal.
\end{lemma}
\begin{proof}
First we show that $M$ is a completely prime ideal of $T$.
Let $f, g\in M$. Since $V$ is a total valuation ring, we can assume that 
$f(0)=g(0)s$ for some $s\in V$. Now if $(f(0)-g(0))u=1$ for some $u\in V$ then \[
g(0)(s-1)u=(g(0)s-g(0))u=(f(0)-g(0))u=1 ,
\] 
which shows that $g(0)\notin M$, a contradiction. 
It is clear that $ft, tf \in M$ for all $f\in M$ and $t\in T$. Thus $M$ is an ideal of $T$.
If $f,g\in M$, then there exists $u, v\in V$ such that $ug(0)=g(0)u=1$ and $vf(0)=f(0)v=1$.
Now we have $(fg)(0)uv=f(0)g(0)uv=f(0)v=1$, so that $fg\notin M$, and $M$ is a completely prime ideal.

Next we show that 
$M$ is a $\tau$-ideal. We note that $J(V)\subseteq M$, and so $M\cap (V-Q)\neq \emptyset$. It is easy to show that $M=J(V)+V[x,\sigma]x$.
Now let $F$ be a finitely generated right $R$-ideal such that $F\subseteq M$. Then there exist $f_1,\dots, f_n\in M-{0}$ with
$F=f_1T+\dots+f_nT$. Since $f_i(0)\in J(V)$ and the set of all right ideal of $V$ is totally ordered, without loss of
generality we can write $f_1(0)V\subseteq f_2(0)V\subseteq\dots \subseteq f_n(0)V$. Thus $f_n(0)$ a right divisor of $f_i(0)$
for all $i$. Now if $f_n(0)\in Q$, then every $s\in J(V)-Q$ is a right divisor of $f_n(0)$. Hence, there exists an 
element $s\in J(V)-Q$ such that $s$ is a right divisor $f_i(0)$ for all $i$. Since $x=ss^{-1}x$, the element $s$ is 
a right divisor of $f_i$ for all $i$. Thus $f_1T+\dots +f_nT\subseteq sT$ and $(f_1T+\dots +f_nT)_{\nu}\subseteq sT\subseteq M$,
which implies that $M$ is a right $\tau$-ideal.

In much the same way one can show that $M$ is a left $\tau$-ideal.
\end{proof}

\begin{remark}\label{IN}
If $*: F_r(S)\rightarrow F_r(S)$ is a star operation on $S$, and if $I_r(S)$ is the set of all integral right 
ideals of $S$, then we have $I\subseteq I^*\subseteq S^*=S$ for all $I\in I_r(S)$. Therefore, each star operation 
on $S$ induces a function $I\rightarrow I^*$ on $I_r(S)$ such that conditions (1)--(3) of Definition~\ref{star} are satisfied.
Moreover, if $F\in F_r(S)$, then $F=q^{-1}I$ for some $I\in I_r(S)$, and $F^*=(q^{-1}I)^*=q^{-1}I^*$. Therefore, a star $*$
on $S$ is completely determined by its action on $I_r(S)$.
\end{remark}

In the following lemma $F^*_r(S)$ denotes the set of all right $*$-ideals of $S$.
\begin{lemma}
A nonempty subset $F'\subset F_r(S)$ satisfies $F'=F^*_r(S)$ for some $*$-operation on $S$ if and only if
all of the following conditions are satisfied:
\begin{itemize}
\item[(1)]
$S\in F'$;
\item[(2)]
$I\in F'$ implies that $uI\in F'$ for each $u\in U(Q)$;
\item[(3)]
$\emptyset \neq \{I_{\alpha}\}\subseteq F'$ with $\cap I_{\alpha}\neq \{0\}$ implies that $\cap I_{\alpha}\in F'$.
\end{itemize}
If (1)--(3) hold we can define $I^*=\cap \{J\in F' : I\subseteq J\}$. 
\end{lemma}
\begin{proof}
Suppose that $*$ is a star operation on $S$. By Definition~\ref{star}
and part (b) of \cite[Proposition~32.2]{G92}, the set $F^*_r(S)$ satisfies (1)--(3). 
We always have $I^*=\cap \{J\in F_r^*(S) : I\subseteq J\}$.
Conversely, for $I\in F_r(S)$ we define $I^*:=\cap \{J\in F' : I\subseteq J\}$.
Since the intersection of any collection of right $S$-ideals is again a right $S$-ideal,
we have $I^*\in F_r(S)$. From (1) and (2) we can conclude that $S^*=S$, and $(uI)^*=uI^*$
for all $I\in F_r(S)$, and $u\in U(Q)$. The definition of $*$ implies that $I\subseteq I^*$
and $I^*\subseteq J^*$, wherever $I\subseteq J$.
Now it is clear from the definition and condition (3) that $F'=\{I\in F_r(S) : I=I^*\}$.
Thus $(I^*)^*=I^*$. 
\end{proof}

\begin{definition}\label{D}
A non-Artinian ring $S$ in a simple Artinian ring $Q$ is called a \textit{discrete} Dubrovin
valuation ring if $S$ is a maximal Dubrovin valuation ring of $Q$ and $J(S)\neq J(S)^2$.  
\end{definition}
A Dubrovin valuation ring $S$ is discrete if and only if $S\cap F$ is a discrete valuation ring of the field $F$,
where $F$ is the center of $Q$, see part c of \cite[Proposition 2.7] {J}. 
There are several equivalent definitions to Definition~\ref{D}, see \cite [Theorem 2.6 ]{J}.

\begin{lemma}
Let $S$ be a discrete Dubrovin valuation ring of a simple Artinian ring $Q$. Then 
$*=i_d$ for all nontrivial star operations $*$ on $S$.

\end{lemma}
\begin{proof}
Since $I\subseteq I^*\subseteq I_{\nu}$, it is enough to prove that every right $S$-ideal is divisorial.
Since $S$ has rank one, every principal right $S$-ideal is a two sided $S$-ideal. That is, 
$aS=Sa$ for all nonzero $a\in Q$ (see the discussion before \cite[Lemma 8]{BMO00}). 
Let $H_r(S)$ and $H(S)$ be the set of all principal right $S$-ideals and principal $S$-ideals respectively. 
Then by part (i) of \cite[Theorem 9]{BMO00}, $H_r(S)=H(S)\cong F_r(S)=F(S)$, 
which show that every right $S$-ideal is divisorial.
\end{proof}
\begin{lemma}
Let $S$ be a Dubrovin valuation ring of $Q$ such that $*=i_d$
for all nontrivial semistar operations $*$ on $S$. Then $S$ is discrete.
\end{lemma}
\begin{proof}
Assume that $J(S)= J(S)^2$. Then, by \cite[Proposition 1.3]{IMU04}, for every $a\in S$ the right ideal 
$aJ(S)$ is not divisorial. Thus $\nu\neq i_d$, which is a contradiction.
Therefore, $J(S)\neq J(S)^2$. Now let $T$ be a proper overring of $S$ in $Q$. We define $I^{*_T}=IT$
for every $S$-submodule $I$ of $Q$. It is easy to see that ${*_T}$ is a semistar on $S$, 
and $I^{*_T}=\cup \{JT : J\subseteq I$, $J$ is finitely generated $\}$.
Thus $*_T$ is semistar and $S^{*_T}=ST=T\neq S$. 
This is a contradiction, because ${*_T}\neq i_d$. Hence $S$ is a maximal Dubrovin valuation ring. 
\end{proof}

\section{Right Pr\"ufer $\nu$-multiplication orders}\label{sec4}

An integral domain in which every finitely generated ideal is invertible is called a Pr\"ufer domain.
An integral domain in which every finitely generated ideal is $t$-invertible 
is called a Pr\"ufer $\nu$-multiplication domain and denoted by P$\nu$MD.
The P$\nu$MDs include Krull domains, Pr\"ufer domains, $\GCD$-domains and unique factorization domains. 
The aim of this section is introduce a non-commutative version of P$\nu$MDs in simple Artinian rings. 
This version includes P$\nu$MD right Bezout orders, right $\GCD$-orders and right Pr\"ufer orders.

For an additive subgroup $I$ of a simple Artinian ring $Q$ we define $O_l(I):=\{q\in Q: qI\subseteq I\}~,
O_r(I):=\{q\in Q: Iq\subseteq I\}$ and $I^{-1}:=\{q\in Q: IqI\subseteq I\}$. 

Recall that an order $S$ in a simple Artinian ring $Q$ is called a right Pr\"ufer ring if any finitely 
generated right $S$-ideal $I$ is left invertible as a right $S$-ideal and right invertible as a left $O_l(I)$-ideal.
More precisely we have:    
\begin{definition}
An order $S$ in a simple Artinian ring $Q$ is a right \textit{Pr\"ufer order} if and only if
$I^{-1}I=S$ and $II^{-1}=O_l(I)$ for any finitely generated right $S$-ideal $I$. 
\end{definition} 
A left Pr\"ufer order is defined similarly. A Pr\"ufer order $S$ is a right as well as a left Pr\"ufer order.
By \cite[Lemma 1.4 and 1.5] {MMU97}, an order $S$ in a simple Artinian ring $Q$ is a right \textit{Pr\"ufer order} if
and only if $(S:I)_lI=S$ and $I(S:I)_l=O_l(I)$ for any finitely generated right $S$-ideal $I$.

\begin{lemma}
Let $S$ be a Pr\"ufer order in a simple Artinian ring $Q$ with finite dimension over its center.
For any nonzero $I\in F_r(S)$ we define $I\rightarrow I^{\omega}=\cap IS_M$, where $M$ runs over all maximal ideals of $S$.
Then:

\begin{itemize}
\item[(i)]
$\omega$ is a star operation and $IS_M=I^{\omega}S_M$. 
\item[(ii)]
If $I$ is a finitely generated right ideal of $S$ such that $(JI)^{\omega}\subseteq (KI)^{\omega}$
and rank $S_M$ is one for every $M$, then $J^{\omega}\subseteq K^{\omega}$.
\item[(iii)]
If $S$ is Noetherian then $\omega=\nu$.

\end{itemize}
\end{lemma}

\begin{proof}

\begin{itemize}
\item[(i)]
By \cite[Theorem 22.8]{MMU97} and \cite[Lemma 2.4]{M93}, for each maximal ideal $M$ of $S$ the localization $S$ at $M$ 
exists and $R_M$ is a Dubrovin valuation ring. Furthermore, $S=\cap S_M$ where $M$ runs over all maximal ideals of $S$.
First we show that $I^{\omega}\in F_r(S)$ for all $I\in F_r(S)$. Then
we prove that $\omega$ satisfies all three conditions (1)--(3) of Definition~\ref{star}.
Since $I^{\omega}S=(\cap _MIS_M)S=\cap_MIS_MS=\cap_MIS_M=I^{\omega}$, the $I^{\omega}$ is a nonzero right $S$-submodule of $Q$.
Now let $u\in U(Q)$ such that $uI\subseteq S$. Then $uI^{\omega}=u(\cap _MIS_M)=\cap _MuIS_M\subseteq \cap _MSS_M=\cap _MS_M=S$.
Therefore, $I^{\omega}\in F_r(S)$.
To prove that $I\rightarrow I^{\omega}$ is a star operation, we need to show that $\omega$ satisfies 
all three conditions (1)--(3) of Definition~\ref{star}.

To prove condition of (1),
let $u\in U(Q)$. Then $u{I^{\omega}}=u\cap_MIS_M=\cap_MuIS_M=\cap_M(uI)S_M=(uI)^{\omega}$, which implies (1).

That condition (2) holds is clear.

To prove condition (3) we can write $I=u^{-1}J$, where $J$ is a right integral $S$-ideal and $u\in U(Q)$. By part (1) we have
$I^{\omega}=u^{-1}J^{\omega}$ and $(I^{\omega})^{\omega}=u^{-1}(J^{\omega})^{\omega}$.
Therefore, we only need to show that $J^{\omega}=(J^{\omega})^{\omega}$.
Since the right ideal $JS_M$ is an extension of the right ideal $J$ of $S$ to $S_M$ and $JS_M\subseteq S_M$,
we have $JS_M=(JS_M\cap S)S_M =(JS_M\cap (\cap_MS_M))S_M\supseteq (\cap_MJS_M)S_M\supseteq {J^{\omega}}S_M$.
 Thus $J^{\omega}=\cap_MJS_M\supseteq \cap_MJ^{\omega}S_M=(J^{\omega})^{\omega}$.
It follows that $I\rightarrow I^{\omega}$ is a ${\omega}$-operation.
From the fact that $JS_M\supseteq J^{\omega}S_M$, we conclude that $IS_M=I^{\omega}S_M$.

\item[(ii)]
From (i) and $(JI)^{\omega}\subseteq (KI)^{\omega}$,
we have $(JI)S_M=(JI)^{\omega}S_M\subseteq (KI)^{\omega}S_M=(KI)S_M$. From the fact that
$S_M$ is a Dubrovin valuation ring, $IS_M$ is a finitely generated right $S_M$-ideal
and using \cite[Corollary 5.5]{MMU97}, we have $IS_M=qS_M$ for some regular element $q\in Q$.
Since $S_M$ is a rank one and $q$ is a unit in $Q$, we have $qS_Mq^{-1}=S_M$.
Thus $(JI)S_M=JqS_M=(JS_M)q\subseteq(KI)S_M=KqS_M=(KS_M)q$,
and so $(JS_M)\subseteq (KS_M)$. Hence $J^{\omega}\subseteq K^{\omega}$. \qedhere
\end{itemize}
\end{proof}

Let $S$ be an order in a simple Artinian ring $Q$ and $*$ is a star operation on $S$.
 A right $*$-ideal $I$ is called of \textit{finite type} if there
exists a finitely generated right $S$-ideal $J$ such that $I=J^*$. Let 
$*=\nu$ and $H_r(R)$ be the set of all right $\nu$-deals 
of finite type.
\begin{definition}

An order $S$ in a simple Artinian ring $Q$ is called a 
\textit{right Pr\"ufer $\nu$-multiplication} order
if \[
((S:I)_lI)^{\tau}=S \quad \text{ and } \quad (I\,(S:I)_l)^{\tau}=O_l(I)
\] for every $I\in H_r(R)$.

A \textit{left Pr\"ufer $\nu$-multiplication} order is defined similarly.
A \textit{Pr\"ufer $\nu$-multiplication} order is simultaneously a right and a left $\nu$-Pr\"ufer
order.
\end{definition}
\begin{lemma}
The following are right Pr\"ufer $\nu$-multiplication orders:
\begin{itemize}
\item[(1)]
Commutative Krull domains;
\item[(2)]
Commutative Pr\"ufer $\nu$-multiplication domains;
\item[(3)]
Right Bezout orders;
\item[(4)]
Right $\nu$-Bezout orders;
\item[(5)]
Right Pr\"ufer orders;

\end{itemize}
\end{lemma}
\begin{proof}
Every commutative Krull domain is a commutative Pr\"ufer $\nu$-multiplication domain. To prove
(1) and (2) it is thus enough to show that the Pr\"ufer $\nu$-multiplication domain $S$ is a right Pr\"ufer $\nu$-multiplication
order. To see this, let $I$ be a right $S$-ideal of finite type. From the fact that $S$ is commutative and  
$I$ has an inverse with respect to $\tau$ multiplication, we have $I[S:I]_l=S=[S:I]_lI$ and $O_l(I)=S$. Thus
$S$ is a right $\nu$-multiplication order in the fraction field $F$ of $S$.

To prove (3) we only need to prove (4), because every right Bezout order is a right $\nu$-Bezout order. To establish (4),
let $S$ be a right $\nu$-Bezout order and $I$ be a right $S$-ideal of finite type. Then $I=J^{\nu}$
for some finitely generated right $S$-ideal $J$. Since $S$ is right $\nu$-Bezout we have $I=J^{\nu}=qS$
for some regular element $q\in Q$. It is easy to show that $[S:I]_l=[S;qS]_l=Sq^{-1}$ and $O_l(I)=qSq^{-1}$.
Therefore, $I^{-1}I=Sq^{-1}qS=S$ and $I^{-1}I=(qS)(Sq^{-1})=qSq^{-1}=O_l(I)$, as desired. 

Proving (5), the set of all finitely generated right $S$-ideals coincides with $H_r(S)$.
Now let $I\in H_r(S)$. Then $[S:I]_lI=S$ and $I[S:I]_l=O_l(I)$. 
Hence $([S:I]_lI)^{\tau}=S$ and $(I[S:I]_l)^{\tau}=O_l(I)$.  
\end{proof}

In \cite[Theorem 1]{D93} Dubrovin proved that every non-commutative 
Pr\"ufer order is a semi-hereditary order. The following example  
taken from \cite{MY2000}, shows that the converse is not necessarily true.
\begin{example}
Let $p$ be an odd prime and $D=\Z_p[t]_{t\Z_p[t]}$. Then $D$ is a local ring with maximal ideal 
$m=tD$. Let $F$ be the quotient field of $\Z_p[t]$ which is also the quotient field of $D$.
We define the automorphism $\sigma$ on $\Z_p[t]$ by $t^{\sigma}=-t$ and $a^{\sigma}=a$ for all $a\in \Z_p$.
Then $\sigma$ can naturally be extended to $D$. Now let $S=D[x]_{xD[x]}$.
 Then $S=\{f(x)g(x)^{-1} : f(x), g(x)\in D[x] \textit ~{with}~ g(0)\neq 0\}$. Define the 
 epimorphism $\Phi$ from $S$ to $F$ by $\Phi (f(x)g(x)^{-1})=f(0)g(0)^{-1}$, and let $R=\Phi^{-1}(D)$. Then 
 $R$ is a valuation ring of rank $2$ and $P_0+mR, P_0$ and $(0)$ are the only prime ideals of $R$.
 The automorphism $\sigma$ is extended to an automorphism of $D[x]$ by 
  $(f(x))^{\sigma}=a_n^{\sigma}x^n+\dots +a_0^{\sigma}$ for any $f(x)=a_nx^n+\dots+a_0\in D[x]$.
  The group $G$ generated by $\sigma$ has order of $2$ and the skew group ring $R\star G$ is a 
  semihereditary order which is not Pr\"ufer.
\end{example}
\begin{question}
Is $R\star G$ a \textit{ right Pr\"ufer $\nu$-multiplication} order?
\end{question}

Recall that a ring $R$ is called a generalised discrete valuation ring, if the set of right ideals are inversely well-ordered
by inclusion. All right ideals of such a ring are two-sided and actually principal right ideal.
Now we give the definition of non-commutative Krull domain in the sense of Brungs.

\begin{definition}
An integral domain $R$ is called a non-commutative Krull domain if there is a family of generalised discrete valuation
domains $V_i~ , i\in I$, satisfying the following:
\begin{itemize}
\item[(1)]
$R=\cap _{i\in I}V_i$ and each $V_i$ is a subring of the same division ring $Q$ such that $Q=Q(V_i)$ for every $i\in I$.
\item[(2)]
Every $a\in R$ is a unit for all but a finite number of $V_i$.
\item[(3)]
Each $V_i$ is satisfies as $R_{P_i}$ for prime ideal $P_i$ of $R$ such that $P_i\cap P_j$ contains no nonzero prime ideal for $i\neq j$. 
\end{itemize} 
\end{definition}
Every one-sided ideal of $R$ is two-sided and $aR=\cap_{i\in I} aV_i$ for all $a\in Q^*$. Thus $F_r(R)=F_l(R)=F(R)$.
We define the operation $b$ by $A^b=\cap AV_i$, where $A\in F_r(R)$. It is easy to show that $b$ is a star operation.
The equivalence relation on $F_r(R)$ can be defined by $A\sim B$ whenever $A^b=B^b$. 
For each $A\in F_r(R)$, we denote by $[A]$ the equivalence class determined by $A$.
Let $D(R)$ be the set of all equivalence classes. Then by \cite[Theorem 2.3]{M75}, $D(R)$ is an abelian group with defined operation
$[A][B]=[A^bB^b]$. $A^b\subseteq A^{\tau}$ for all $A\in F_r(R)$. 
Hence non-commutative Krull domains are examples of non-commutative Pr\"ufer $\nu$-multiplication domains.
\begin{example}
Let $F$ be a field of characteristic zero and $A:=F[x,y]$ with $xy-yx=1$. Then $A[z]$ is an
Asano order in $Q(A[z])$ \cite{HL72}. For every non unit $a\in A$ the right ideal $I:=aA[z]+zA[z]$ generated by $a, z$
is not projective. Thus $A[z]$ is not a right semi hereditary order and hence not a right Pr\"ufer order. On the 
other hand, $A[z]$ is a Krull domain in the sense of Brungs and hence a Pr\"ufer $\nu$-multiplication order. 
\end{example}

Recall that two orders $R$ and $S$ in a simple Artinian ring $Q$ is called equivalent if there
exist regular elements $a_1, a_2, b_1, b_2\in Q$ such that $a_1Rb_1\subseteq S$ and $a_2Sb_2\subseteq R$.   
\begin{lemma}
Let $S$ be a bounded Krull ring in the sense of Marubayashi. Then $O_l(I)$ is a bounded
Krull ring for any divisorial right $S$-ideal $I$. 
\end{lemma}
\begin{proof}
By \cite[Lemma 2.5]{M76} $O_l(I)$ is a maximal order equivalent to $S$, and by \cite[Theorem 2.6]{M76} 
$O_l(I)$ is a bounder Krull ring.
\end{proof}

\begin{lemma}
Let $S$ be a Noetherian bounded Krull ring in the sense of Marubayashi. Then $S$ is a right Pr\"ufer 
$\nu$-multiplication order.
\end{lemma} 
\begin{proof}
By \cite[Corollary 1.4] {M76} $S$ is a maximal order in the sense of Asano. Now let $I$ be right $S$-ideal of finite
type. Then, by part 4 of \cite[Lemma2] {M78}, $((S:I)_lI^{\nu})^{\nu}=S$ and $(I^{\nu}(S:I)_l)^{\nu}=O_l(I)$.
The ring $S$ is Noetherian and $O_l(I)$ is equivalent to $S$. Hence $O_l(I)$ is also Noetherian. Therefore, 
$(S:I)_lI^{\nu}$ and $I^{\nu}(S:I)_l$ are finitely generated right ideals of $S$ and $O_l(I)$ respectively.
Thus $((S:I)_lI)^{\tau}=((S:I)_lI^{\nu})^{\nu}=S$ and $(I(S:I)_l)^{\tau}=(I^{\nu}(S:I)_l)^{\nu}=O_l(I)$, as desired.
\end{proof}

The following is a generalisation of \cite[Proposition 2.6]{MMU97} to right Pr\"ufer $\nu$-multiplication orders.
\begin{proposition}
Let $S$ be a right Pr\"ufer $\nu$-multiplication order in simple Artinian ring $Q$ and 
$T$ is an overring of $S$. Then $T$ is also a Pr\"ufer $\nu$-multiplication order in $Q$.
\end{proposition}

\begin{proof}
Let $I'\in H_r(T)$. Then $I'=(a_1T+\cdots+a_nT)^{\nu}$ for some finitely generated right $T$-ideal of $a_1T+\cdots+a_nT$.
Put $I=(a_1S+\cdots+a_nS)^{\nu}$. Then $I\in H_r(S)$ and $((S:I)_lI)^{\tau}=S$ and $(I(S:I)_l)^{\tau}=O_l(I)$.
Now let $x\in (S:I)_l$. Then $x(a_1S+\cdots+a_nS)\subseteq xI\subseteq S$ and so $x(a_1S+\cdots+a_nS)T\subseteq xIT\subseteq ST$.
Thus $x(a_1T+\cdots+a_nT)\subseteq T$ and $xI'\subseteq T$, which shows that $(S:I)_l\subseteq (T:I')_l$. From
$(S:I)I\subseteq (T:I')_lI'\subseteq T,~I(S:I)\subseteq I'(T:I')_l\subseteq O_l(I')$
and since $S$ is a right Pr\"ufer $\nu$-multiplication order, we have
$S=((S:I)I)^{\tau}\subseteq ((T:I')_lI')^{\tau}\subseteq T$ and $O_l(I)=(I(S:I))^{\tau}\subseteq (I'(T:I')_l)^{\tau}\subseteq O_l(I')$.
 Hence $1\in ((T:I')_lI')^{\tau}$ and $1\in (I'(T:I')_l)^{\tau}$. Therefore, $((T:I')_lI')^{\tau}=T$ and $(I'(T:I')_l)^{\tau}=O_l(I')$.  
\end{proof}

\begin{remark}
By \cite[Lemma 1.5]{MMU97} a right $S$-ideal $I$ of $Q$ is a projective $S$-module if and only if 
$I(S:I)_l=O_l(I)$. Thus if $I$ is projective $S$-module, we always have $(I(S:I)_l)^{\tau}=O_l(I)$ but the converse 
is not true. For example, the Krull ring $S=\F[x,y]$ is not a semi-hereditary.
Thus there exists a finitely generated $S$-ideal $I$ which is not projective. 
Now, since $S$ is a Pr\"ufer $\nu$-multiplication domain, 
we have $(I(S:I)_l)^{\tau}=O_l(I)$. In \cite[Proposition 2.5]{MMU97}
it is proven that a right Pr\"ufer order is left Pr\"ufer and vice versa. The proof relies on the
semi-hereditarity of Pr\"ufer orders and \cite[Lemma 1.5]{MMU97}. 
The Pr\"ufer $\nu$-multiplication order is not necessarily semi-hereditary. 
As yet we have been unable to prove or disprove that a 
right Pr\"ufer $\nu$-multiplication order is a left
$\nu$-multiplication order. 
\end{remark}

\bibliographystyle{amsplain}

\end{document}